\documentclass[letterpaper, 10 pt, conference]{ieeeconf}  

\IEEEoverridecommandlockouts                              
\usepackage{amsmath, amssymb, amsxtra, mathrsfs}
\usepackage{graphics, epsfig, here, epstopdf, fancyhdr,graphicx}
\usepackage{bbm}
\usepackage{color}
\usepackage{multirow}
\usepackage[all]{xy}

\newcommand{\cP}{\mathcal{P}}
\newcommand{\cC}{\mathcal{C}}

\newcommand{\cJ}{\mathcal{J}}

\newcommand{\cS}{\mathcal{S}}

\newcommand{\fA}{\,\forall\,}

\newcommand{\dS}{\displaystyle}

\newcommand{\bq}[1]{{\left[#1\right]}}

\newcommand{\cG}{\mathcal{G}}

\newcommand{\thalf}{\tfrac{1}{2}}

\newcommand{\rd}{\mathrm{d}}

\usepackage{algorithmic}
\usepackage{epstopdf}
\usepackage{color}
\usepackage{dashrule}
\usepackage{booktabs}

\title{\LARGE \bf
Efficient Dynamic Compressor Optimization \\ in Natural Gas Transmission Systems
}

\author{Terrence W. K. Mak, Pascal Van Hentenryck, Anatoly Zlotnik, Hassan Hijazi, and Russell Bent
\thanks{T. W. K. Mak and H. Hijazi are at NICTA and Australian National University, Canberra, ACT, Australia. Email: \{Terrence.Mak $\mid$ Hassan.Hijazi\}@nicta.com.au;
P. Van Hentenryck is at the Department of Industrial \& Systems Engineering, University of Michigan, Ann Arbor, MI, USA. Email: pvanhent@umich.edu; R. Bent and A. Zlotnik are at Los Alamos National Laboratory, Los Alamos, NM, USA. Email: \{rbent $\mid$ azlotnik\}@lanl.gov.
}
}

\begin{document}

\maketitle
\thispagestyle{empty}
\pagestyle{empty}

\begin{abstract}
The growing reliance of electric power systems on gas-fired generation to balance intermittent sources of renewable energy has increased the variation and volume of flows through natural gas transmission pipelines.  Adapting pipeline operations to maintain efficiency and security under these new conditions requires optimization methods that account for transients and that can quickly compute solutions in reaction to generator re-dispatch. This paper presents an efficient scheme to minimize compression costs under dynamic conditions where deliveries to customers are described by time-dependent mass flow. The optimization scheme relies on a compact representation of gas flow physics, a trapezoidal discretization in time and space, and a two-stage approach to minimize energy costs and maximize smoothness. The resulting large-scale nonlinear programs are solved using a modern interior-point method. The proposed optimization scheme is validated against an integration of dynamic equations with adaptive time-stepping, as well as a recently proposed state-of-the-art optimal control method.  The comparison shows that the solutions are feasible for the continuous problem and also practical from an operational standpoint.  The results also indicate that our scheme provides at least an order of magnitude reduction in computation time relative to the state-of-the-art and scales to large gas transmission networks with more than 6000 kilometers of total pipeline.
\end{abstract}

\section{Introduction}
\label{sec:intro}

In the last two decades, the increasing integration of renewable
sources of energy in the electrical power system, together with the
recent availability of natural gas in the United States, acted as a
primary driver of installation of gas-fired electric power plants to
meet most of the demand for new generating capacity \cite{lyons13}.
Gas-fired generators may go online and shut down several times a day, and
can rapidly modify their outputs, making them attractive generation
resources to balance the fluctuations of renewable energy sources such
as wind and solar \cite{li08,mitei14}.

Historically, withdrawals from natural gas transmission systems came
from utilities and industrial consumers, which are highly predictable
and exhibit low variation in demand \cite{chertkov14}. These withdrawals
are traded using day-ahead contracts for fixed deliveries and
implicilty assume that injections and withdrawals remain nearly
constant. As a result, natural gas transmission systems could be
largely controlled and optimized using steady-state modeling
\cite{luongo91,riosmercado15}. Early studies
\cite{wong68,rothfarb70,luongo91} focused on optimizing steady-state
gas flows, for which the state equations are algebraic relations.
Recent efforts have improved and scaled up optimization techniques for
similar problems
\cite{midthun09,borraz10phd,misra15,babonneau14,borraz15c}.  In
short-term operations, the operating set-points for gas compressor
stations can be readily changed, and compressor optimization for
steady-state flows has been solved in the form of an optimal gas flow
(OGF) \cite{misra15}.

However, it is no longer appropriate to restrict attention to
steady-state approximations, which cannot adequately describe the
physics of high volume gas flows that may fluctuate significantly
throughout the day according to gas-fired generator dispatch and
commitment schedules \cite{peters12a,peters12b}. Recent studies have
highlighted the challenges created by such increasing flow variability and
volume \cite{chertkov15hicss}, which are predominantly caused by
the growing use of gas-fired power plants for electricity
generation \cite{lyons13,levitan14}. The integration issues that face
the increasingly interdependent electric and gas systems have created
growing concerns in both industry sectors \cite{mitei14}.  As a
result, in order to enable natural gas systems to inter-operate with
electric power systems on the time-scale of generator dispatch, the
OGF must be extended to account for transient flow conditions.  {\em
  Therefore, new optimization techniques that use dynamic models of
  gas flows on pipeline networks are required} \cite{shahidehpour14}.

An automatic control methodology for optimally managing transient
flows in gas transmission systems will require a stable, accurate,
physics-based, and rapid technique for computing model-based
compressor control protocols.  Automatic grid control tools are
already industry standards for management of generator dispatch for
electric power systems \cite{pjm15}.  Such power system optimization
models consider one or a sequence of steady states
\cite{litvinov10}\footnote{Recently, efficient optimization methods
  for transient repair and restoration of power systems
  \cite{hijazi15,vanhentenryck15} have been facilitated by
  discretizations characterized by sparse constraints and appropriate
  nonlinear relaxations.}, but the significant qualitative difference
in the physics of gas pipeline networks prevents these techniques from
being directly transferred. Indeed, the dynamics in the electricity
and gas networks operate at fundamentally different space and time
scales.

Gas pipeline flow dynamics on the relevant spatial and temporal scales
do not experience waves or shocks, and can be represented by the Euler
equations for compressible gas flow in one-dimension with significant
simplifications \cite{osiadacz84,herty10}.  These partial differential
equations (PDEs) are highly nonlinear, however, and are challenging to
simulate \cite{thorley87}, particularly when instances over many
domains are coupled at the boundaries over a network.  This
nonlinearity and complexity of gas pipeline network dynamics is also
an obstacle to tractable optimization of such flows under transient
conditions.  Several studies have proposed optimization schemes for
gas networks on the time-scale of daily operations
\cite{rachford00,ehrhardt05,steinbach07pde,zlotnik15cdc}, {\em and the
  issues of computation time and scalability have been noted
  repeatedly.}

In this manuscript, we consider the Dynamic Optimal Gas Flow (DOGF) problem, which
generalizes the OGF to capture the dynamics of a gas pipeline network. The
objective of the DOGF is to minimize the cost of gas compression required for
pipeline operations subject to system pressure constraints, and whose deliveries to customers are described by time-dependent mass flow. {\em Our main contribution is an efficient optimization scheme for the DOGF, which is validated
by an accurate simulation method for gas pipeline networks
with dynamic flows and compressors.}

The key aspects of our optimization scheme can be summarized as
follows. First, following \cite{zlotnik15cdc}, the hydrodynamic
relations that describe gas flow are discretized in time and space
using first order approximations, and several relaxations of the
resulting nonlinear constraints are employed. Second, in contrast to
earlier work that was based on pseudospectral discretization, our proposed scheme uses a trapezoidal discretization but compensates
for the potential loss in accuracy by a two-stage optimization
approach. In the first stage, our scheme optimizes the compression
cost (the original objective).  In order to obtain a physically realistic solution, the second stage optimizes the smoothness of the compressor ratios while ensuring that the overall compression costs remain close to the value found in the first stage.
The resulting large-scale nonlinear optimization problems (with up to more than 130,000 decision variables) are solved using the state-of-the-art IPOPT system \cite{wachter06on}, enabling us to exploit the significant progress in nonlinear optimization in
the last decades.

The solutions produced by our optimization scheme are verified in two
ways:
\begin{enumerate}
\item By comparing them to a validated dynamic simulation method for
  gas pipeline networks with transient compression
  \cite{grundel13b,zlotnik15dscc}, which is seeded with the compressor
  ratios of our optimized solutions; and

\item By comparing them to the state-of-the-art method for the DOGF.
\end{enumerate}
The validation process indicates that our optimization scheme produces
solutions with no pressure violations and with physically meaningful mass flow
and pressure trajectories in the corresponding simulation. Moreover, the
compressor ratios in our solutions only have negligible differences
compared to the state-of-the-art method. The main benefit of our
optimization scheme, however, is its computational efficiency: Compared
to earlier approaches, {\em it reduces solution times by one to two
  orders of magnitude}. In particular, it solves a previously
investigated 24-pipe gas network case study in less than 30 seconds
and demonstrates scalability on networks with 24, 40 and 135 pipes,
with total pipeline lengths of 477, 1118, and 6964 kilometers, respectively.

The rest of the paper is organized as follows. Section
\ref{sec:problem} summarizes the physical modeling of gas pipeline
networks and poses the DOGF.  Section \ref{sec:scheme} describes the
discretization scheme and relaxations required for an efficient
optimization formulation. Section \ref{sec:results} presents
our two-stage optimization approach and implementation details,
followed by computational results for the three case studies, including convergence experiments and validating simulations. Section \ref{sec:conc} concludes with a discussion and future extensions.

\section{Compressor Optimization in Gas Pipelines} \label{sec:problem}

A gas pipeline network can be represented as a directed graph $\cG=(\cJ,\cP)$, where edges $\{i,j\}\in \cP$ represent pipes $P_{ij}$ that connect nodes $i,j\in\cJ$ that represent joints $J_i$ and $J_j$.  The dynamic state on the pipe $P_{ij}$ is given by pressure $p_{ij}$ and mass flow $q_{ij}$, which evolve on a time interval $[0, T]$ and the distance variable $x_{ij} \in[0, L_{ij}]$, where $L_{ij}$ is the length of pipe $P_{ij}$.  We are interested in the subsonic and isothermal regime of transients that do not excite shocks or waves, i.e., where the flow velocity through a pipe is less than the speed of sound $a$ in the gas, and temperature is assumed to be constant.  The flow dynamics on a single pipe $P_{ij}$ can be adequately described in this regime \cite{herty10} by
\begin{align}
\frac{\partial p_{ij}}{\partial t} + \frac{a^2}{A_{ij}} \frac{\partial q_{ij}}{\partial x} = 0 \label{eq:orig1}\\
2 p_{ij} \frac{\partial p_{ij}}{\partial x} + \frac{\lambda a^2}{D_{ij} A_{ij}^2} q_{ij} \lvert q_{ij} \rvert = 0 \label{eq:orig2}
\end{align}
The parameters for each pipe $P_{ij}$ are the pipe diameter $D_{ij}$ and cross-sectional area $A_{ij}$, and the speed of sound $a$ and friction factor $\lambda$ are assumed uniform and constant throughout the system.  The second term in \eqref{eq:orig2} approximates friction effects.  The gas dynamics on a pipeline segment are represented using \eqref{eq:orig1}-\eqref{eq:orig2} and possess a unique solution when any two of the boundary conditions $p_{ij}(t,0)$, $q_{ij}(t,0)$, $p_{ij}(t,L_{ij})$, or $q_{ij}(t,L_{ij})$ are specified.  For both computational and notational purposes, we apply a transformation to dimensionless variables \cite{herty10} given by
\begin{align}
\tilde{p}_{ij} &= \frac{p_{ij}}{p_N}, &\quad \tilde{q}_{ij} &= \frac{q_{ij}}{q_N}, \label{dimensionless} \\
\tilde{x}_{ij} &= x \frac{\lambda a^2 q_N^2}{D_{ij}A_{ij}^2p_N^2}, & \quad
\tilde{t}_{ij} &= t \frac{\lambda a^4 q_N^3}{D_{ij}A_{ij}^3p_N^3}, \nonumber
\end{align}
where $p_N$ and $q_N$ are scaling constants.  This results in the dimensionless equations
\begin{align}
\frac{\partial p_{ij}}{\partial t_{ij}} + \frac{\partial q_{ij}}{\partial x_{ij}} = 0, \label{eq:norm1}\\
2 p_{ij} \frac{\partial p_{ij}}{\partial x_{ij}} + q_{ij} \lvert q_{ij} \rvert = 0, \label{eq:norm2}
\end{align}
where we omit the $\sim$ label over the variables for readability.  Note that the space and time variables $x_{ij}$ and $t_{ij}$ are now pipe-dependent.  We write these variables as $x$ and $t$ when it is clear from the context that they correspond to a specific pipe $P_{ij}$.  Design limits and regulations for pipeline systems require pressure to remain within specified bounds given by
\begin{align} \label{presbounds}
\underline{p}_{ij} \leq p_{ij}(t,x) \leq \overline{p}_{ij}.
\end{align}


The momentum dissipation due to the friction term in \eqref{eq:orig2} causes the gas pressure to decrease, hence it must be augmented by compressors to maintain the minimum required pressure.  We define $\cC\subset\cP$ as the subset of pipes that have compressors.  The action of compressors is modeled as conservation of flow and an increase in pressure at a point $c_{ij}\in[0, L_{ij}]$ by a multiplicative ratio $R_{ij}(t)>0$ that may depend on time.  Specifically,
\begin{align}
\lim_{x\searrow c_{ij}}p_{ij}(t,x)&=R_{ij}(t)\lim_{x\nearrow c_{ij}}p_{ij}(t,x), \label{eq:compbal1}\\
\lim_{x\searrow c_{ij}}q_{ij}(t,x)&=\lim_{x\nearrow c_{ij}}q_{ij}(t,x). \label{eq:compbal2}
\end{align}
The cost of compression $S_{ij}$ is proportional to the required power \cite{misra15}, and is approximated by
\begin{align} \label{obj1}
S_{ij}(t) = \eta^{-1}|q_{ij}(t,c_{ij})|(\max\{R_{ij}(t),1\}^{2K}-1)
\end{align}
with $0<K=(\gamma-1)/\gamma<1$, where $\gamma$ is the heat capacity ratio and $\eta$ is a compressor efficiency factor.  In this study we do not consider pressure regulation (decompression), so the compressor ratio for a given station must remain bounded within a feasible operating region
\begin{align} \label{compconst}
\max\{\underline{R}_{ij},1\} \leq R_{ij} \leq \overline{R}_{ij}.
\end{align}

In addition to the dynamic equations \eqref{eq:norm1}-\eqref{eq:norm2} and continuity conditions for compressors \eqref{eq:compbal1}-\eqref{eq:compbal2} that characterize the system behavior on each pipe $P_{ij}\in\cP$, we specify balance conditions for each joint $J_i\in\cJ$.  We first define variables for the unique nodal pressure $p_i(t)$ at each junction, as well as mass flow injections $d_i(t)$ from outside the system (negative for consumptions/withdrawals). Each joint $J_j\in \cJ$ then has a flow balance condition
\begin{align} \label{bal1}
\!\!\!\! \sum_{J_i\in \cJ: P_{ij}\in \cP} \!\! q_{ij}(t,L_{ij}) - \!\!\!\! \sum_{J_k\in \cJ: P_{jk}\in \cP} \!\! q_{jk}(t,0) = f_j(t),
\end{align}
as well as a pressure continuity condition
\begin{align} \label{bal2}
p_{ij}(t,L)=p_j(t)=p_{jk}(t,0), \qquad \qquad  \\ \qquad \qquad \fA J_i,J_k\in \cJ \text{ s.t. } P_{ij},P_{jk}\in\cP. \nonumber
\end{align}
A subset of the junctions $\cS\subset \cJ$ may be treated as ``slack'' nodes, which reasonably represent large sources of gas to a transmission system.  For these junctions, the mass inflow $f_i(t)$ is a free variable and the nodal pressure is a parameter
\begin{align} \label{bal3}
p_i(t)=s_i(t).
\end{align}
For the remaining junctions, which reasonably represent consumers or small suppliers of gas to the system, the nodal pressure $p_i(t)$ free and the mass inflow is a parameter
\begin{align} \label{bal4}
f_i(t)=d_i(t).
\end{align}

The optimization problem that we aim to solve involves a gas pipeline network for which the conditions at each joint are parameterized by an injection/withdrawal $d_i(t)$ or supply pressure $s_i(t)$.  The design goal is for the system to deliver all of the required flows $d_i(t)$ while maintaining feasible system pressure given the physics-based dynamic constraints, and the objective is to minimize the cost of compression over a time interval $[0,T]$.  This cost objective is given by
\begin{align}
C=\sum_{P_{ij}\in\cC} \int_0^T S_{ij}(t) \rd t. \label{cost1}
\end{align}
In this study we consider time-periodic boundary conditions on the system state and controls, i.e.,
\begin{align}
&\!\!\!\!\!\!\!  p_{ij}(0,x)=p_{ij}(T,x),  q_{ij}(0,x)=q_{ij}(T,x),  \fA P_{ij}\in\cP  \!\label{periodicstate}\\
&\!\!\!\!\!\!\!  R_{ij}(0)=R_{ij}(T), \quad \fA P_{ij} \in \cC \label{periodiccontrol}
\end{align}
and therefore feasible parameter functions also must satisfy $d_i(0)=d_i(T)$ and $s_i(0)=s_i(T)$.  The formulation is
\begin{equation} \label{prob:ocp}
\begin{array}{llll}
\min & C \text{ in } \eqref{cost1} \\ s.t.
& \text{pipe dynamics: }  \eqref{eq:norm1},\eqref{eq:norm2}  \\
& \text{compressor continuity: } \eqref{eq:compbal1},\eqref{eq:compbal2} \\
& \text{joint conditions: }  \eqref{bal1},\eqref{bal2}\\
& \text{density \& compression constraints: } \eqref{presbounds},\eqref{compconst} \\
& \text{periodicity constraints: } \eqref{periodicstate},\eqref{periodiccontrol} \\
& \text{injection parameters: } \eqref{bal3},\eqref{bal4} \\
& \text{compressor power: } \eqref{obj1}
\end{array}
\end{equation}

In the next section, we describe a spatial and temporal discretization scheme and relaxation conditions that facilitate efficient solution of this PDE-constrained optimization problem using standard nonlinear programming tools.

\section{Discretization to a Nonlinear Program} \label{sec:scheme}

We create a discretization scheme to balance the high nonlinearity in the spatiotemporal dynamics \eqref{eq:norm1}-\eqref{eq:norm2} among a collection of auxiliary variables in which the constraints in Problem \eqref{prob:ocp} possess a sparse representation.

For each pipe $P_{ij}\in\cP-\cC$, we create a set of $M+1$ time points $t_m^{ij}$ and $N_{ij}+1$ space points $x_n^{ij}$ defined by
\begin{align}
&t_m^{ij}=m\Delta_{ij}^t, \quad m=0,1,\ldots,M, \label{tdisc} \\
&x_n^{ij}=n\Delta_{ij}^x, \quad m=0,1,\ldots,N_{ij}, \label{xdisc} \\
&\Delta_{ij}^t=\frac{T_{ij}}{M}, \quad \Delta_{ij}^x=\frac{L_{ij}}{N_{ij}}. \label{discsteps}
\end{align}
Here $\Delta_{ij}^t$ and $\Delta_{ij}^x$ are (dimensionless) time and space discretization steps, and $T_{ij}$ is the dimensionless time horizon for pipe $P_{ij}$ obtained from $T$ according to \eqref{dimensionless}.  We omit the subscripts $\{ij\}$ on $N_{ij}$ when they are clear from the context.  For each of $(M+1)\times(N_{ij}+1)$ discrete points in $\{(t_m,x_n):0\leq m \leq M,0\leq n \leq N_{ij}\}$ within the (dimensionless) domain $[0,T_{ij}]\times[0,L_{ij}]$ for the flow dynamics on a pipe $P_{ij}$, we define
\begin{align}
p_{ij}^{mn} \approx p_{ij}(t_m,x_n), \quad q_{ij}^{mn} \approx q_{ij}(t_m,x_n) \label{varsdef}
\end{align}
to be the pressure and mass flow variables at time $t=t_m$ and location $x=x_n$.  We also define temporal and spatial derivative variables at time $t=t_m$ and location $x=x_n$ by
\begin{align}
&pt_{ij}^{mn} \approx \frac{\partial p_{ij}}{\partial t}(t_m,x_n), \quad
px_{ij}^{mn} \approx \frac{\partial p_{ij}}{\partial x}(t_m,x_n), \label{varsderiv1} \\
& \qquad \qquad \qquad qx_{ij}^{mn} \approx \frac{\partial q_{ij}}{\partial x}(t_m,x_n). \label{varsderiv2}
\end{align}
A constraint that relates the discrete variables \eqref{varsdef} to their derivatives \eqref{varsderiv1}-\eqref{varsderiv2} is created by approximating the integral over a time or space step by the trapezoid rule.  This yields
\begin{align}
\forall P_{ij} \in \cP-\cC, 0 \leq m \leq M-1, 0 \leq n \leq N: \nonumber \\
p_{ij}^{m+1,n} - p_{ij}^{mn} \approx \frac{\Delta_{ij}^t}{2} (pt_{ij}^{m+1,n} + pt_{ij}^{mn} )  \label{eq:Trape1} \\
\forall P_{ij} \in \cP-\cC, 0 \leq m \leq M, 0 \leq n \leq N-1: \nonumber \\
p_{ij}^{m,n+1} - p_{ij}^{mn} \approx \frac{\Delta_{ij}^x}{2} (px_{ij}^{m,n+1} + px_{ij}^{mn} )  \label{eq:Trape2} \\
q_{ij}^{m,n+1} - q_{ij}^{mn} \approx \frac{\Delta_{ij}^x}{2} (qx_{ij}^{m,n+1} + qx_{ij}^{mn} )  \label{eq:Trape3}
\end{align}
The nondimensional dynamic equations \eqref{eq:norm1}-\eqref{eq:norm2} are then discretized in the above variables according to
\begin{align}
\forall P_{ij} \in P - C, 0 \leq m \leq M,& 0 \leq n \leq N: \nonumber \\
pt_{ij}^{mn} + qx_{ij}^{mn} &= 0 \label{eq:Dis1}\\
2 p_{ij}^{mn} px_{ij}^{mn} + q_{ij}^{mn} \lvert q_{ij}^{mn} \rvert &= 0 \label{eq:Dis2}
\end{align}

For each pipe $P_{ij}\in\cC$, we define the discrete compression variables $R_{ij}^m$ for $m=0,1,\ldots,M$, and suppose that the compressor is located at $c_{ij}=x_k$ for some $0\leq k \leq N$, where the dependence of $k$ on the pipe $P_{ij}$ in question is clear from the context.  The pipe is then divided into two pipes $P_{iju}$ and $P_{ijl}$, with non-dimensional lengths $L_{iju}$ and $L_{ijl}$, and for which we we define discrete variables
\begin{align}
p_{iju}^{mn} \approx p_{ij}(t_m,x_n), \, q_{iju}^{mn} \approx q_{ij}(t_m,x_n), \,\, 0\leq n \leq k \label{varsdef2a} \\
p_{ijl}^{mn} \approx p_{ij}(t_m,x_n), \, q_{ijl}^{mn} \approx q_{ij}(t_m,x_n), \,\,\, k\leq n \leq N \label{varsdef2a}
\end{align}
and corresponding spatial derivative variables $pt_{iju}^{mn}$, $px_{iju}^{mn}$, and $qx_{iju}^{mn}$ for $0\leq n\leq k$ and $pt_{ijl}^{mn}$, $px_{ijl}^{mn}$, and $qx_{ijl}^{mn}$ for $k\leq n\leq N$.  These state and derivative variables satisfy
\begin{align}
&\!\!\!\! p_{iju}^{m+1,n} - p_{iju}^{mn} \approx \frac{\Delta_{ij}^t}{2} (pt_{iju}^{m+1,n} + pt_{iju}^{mn} ), \,\, 0\leq n \leq k,  \label{eq:Trape4} \\
&\!\!\!\! p_{ijl}^{m+1,n} - p_{ijl}^{mn} \approx \frac{\Delta_{ij}^t}{2} (pt_{ijl}^{m+1,n} + pt_{ijl}^{mn} ),  \, k\leq n \leq N \label{eq:Trape5}
\end{align}
for $P_{ij} \in \cC$ and $0 \leq m \leq M-1$, and
\begin{align}
&\!\!\!\! p_{iju}^{m,n+1} - p_{iju}^{mn} \approx \frac{\Delta_{ij}^x}{2} (px_{iju}^{m,n+1} + px_{iju}^{mn} ),  \,\, 0\leq n < k,  \label{eq:Trape6} \\
&\!\!\!\! p_{ijl}^{m,n+1} - p_{ijl}^{mn} \approx \frac{\Delta_{ij}^x}{2} (px_{ijl}^{m,n+1} + px_{ijl}^{mn} ),  \, k\leq n \leq N,  \label{eq:Trape7} \\
&\!\!\!\! q_{iju}^{m,n+1} - q_{iju}^{mn} \approx \frac{\Delta_{ij}^x}{2} (qx_{iju}^{m,n+1} + qx_{iju}^{mn} ),  \,\, 0\leq n < k,  \label{eq:Trape8} \\
&\!\!\!\! q_{ijl}^{m,n+1} - q_{ijl}^{mn} \approx \frac{\Delta_{ij}^x}{2} (qx_{ijl}^{m,n+1} + qx_{ijl}^{mn} ),  \, k\leq n \leq N  \label{eq:Trape9}
\end{align}
for $P_{ij} \in \cC$ and $0 \leq m \leq M$.  In addition, we require continuity constraints at the compressor location to connect pipes $P_{iju}$ and $P_{ijl}$, which take the form
\begin{align}
R_{ij}^m=\frac{p_{ijl}^{mk}}{p_{iju}^{mk}}, \,\,
q_{ijl}^{mk}=q_{iju}^{mk},  \label{eq:compcont}
\end{align}
for all $P_{ij} \in \cC$ and $0 \leq m \leq M$. The equations \eqref{eq:norm1}-\eqref{eq:norm2} on either side of the compressor are discretized for $P_{ij}\in\cC$ by
\begin{align}
pt_{iju}^{mn} + qx_{iju}^{mn} &= 0, \,\, 0 \leq n \leq k, \label{eq:Dis3}\\
2 p_{iju}^{mn} px_{iju}^{mn} + q_{iju}^{mn} \lvert q_{iju}^{mn} \rvert &= 0, \,\, 0 \leq n \leq k, \label{eq:Dis4} \\
pt_{ijl}^{mn} + qx_{ijl}^{mn} &= 0, \,\, k \leq n \leq N\label{eq:Dis5}\\
2 p_{ijl}^{mn} px_{ijl}^{mn} + q_{ijl}^{mn} \lvert q_{ijl}^{mn} \rvert &= 0, \,\, k \leq n \leq N \label{eq:Dis6}
\end{align}
for all $P_{ij} \in \cC$ and $0 \leq m \leq M$.  The equations \eqref{eq:Trape1}-\eqref{eq:Dis2} and \eqref{eq:Trape4}-\eqref{eq:Dis6} discretize the dynamic equations \eqref{eq:norm1}-\eqref{eq:norm2} and continuity conditions for compressors \eqref{eq:compbal1}-\eqref{eq:compbal2}.  The pressure variables must also satisfy for all $0\leq m \leq M$ the constraints
\begin{align}
\underline{p}_{ij} \leq p_{ij}^{nm} \leq \overline{p}_{ij},&& P_{ij} \in P - C, 0 \leq n \leq N, \label{eq:PBound1}\\
\underline{p}_{ij} \leq p_{iju}^{nm} \leq \overline{p}_{ij},&& P_{ij} \in C, 0 \leq n \leq k, \label{eq:PBound2}\\
\underline{p}_{ij} \leq p_{ijl}^{nm} \leq \overline{p}_{ij},&& P_{ij} \in C, k \leq n \leq N \label{eq:PBound3}
\end{align}
which discretize \eqref{presbounds}.  In addition, the compression ratio must satisfy
\begin{align} \label{compdisc}
\max\{\underline{R}_{ij},1\} \leq R_{ij}^{m} \leq \overline{R}_{ij}.
\end{align}
for $P_{ij}\in \cC$ and $0\leq m \leq M$.  The cost of compression is then expressed by a constraint
\begin{align} \label{obj2}
S_{ij}^m = \eta^{-1}qm_{ij}^{m}((R_{ij}^m)^{2K}-1)
\end{align}
for $P_{ij}\in \cC$ and $0\leq m \leq M$, where $qm_{ij}^{m}$ is an auxiliary variable with the constraints
\begin{align} \label{qabsconst}
qm_{ij}^{m} \geq q_{iju}^{mk}, \quad  qm_{ij}^{m} \geq - q_{iju}^{mk},
\end{align}
so that $qm_{ij}^{m} = |q_{iju}^{mk}|$ when $R_{ij}^m>1$.  Compressor cost is also constrained to be positive, i.e.,
\begin{align} \label{poscost}
S_{ij}^m \geq 0.
\end{align}

The balance conditions at joints are given in the discrete variables for all $0\leq m \leq M$ and $J_j\in\cJ$ by
\begin{align}
\sum_{J_i\in \cJ: P_{ij}\in \cP} \!\! q_{ij}^{mN} - \!\!\!\! \sum_{J_k\in \cJ: P_{jk}\in \cP} \!\! q_{jk}^{m0} \qquad \qquad \label{bal5}\\ \qquad + \sum_{J_i\in \cJ: P_{ij}\in \cP} \!\! q_{ijl}^{mN} - \!\!\!\! \sum_{J_k\in \cJ: P_{jk}\in \cP} \!\! q_{jku}^{m0}  = f_j^m, \nonumber
\end{align}
and
\begin{align} \label{bal6}
p_{ij}^{mN}=p_j^m=p_{jk}^{m0}, \,\, \fA J_i,J_k\in \cJ \text{ s.t. } P_{ij},P_{jk}\in\cP.
\end{align}
Parametrization of these balance conditions for $0\leq m\leq M$ is given by
\begin{align}
&f_i^m=d_i(t_m), &\quad J_i\in\cJ-\cS, \label{bal7}\\
&p_i^m=s_i(t_m), &\quad J_i\in\cS, \label{bal8}
\end{align}
where $d_i(t)$ and $s_i(t)$ are given flow injection or supply pressure functions.  The time-periodic boundary conditions on the states and controls are given for $0\leq m \leq M$ by
\begin{align}
&p_{ij}^{0n} = p_{ij}^{Mn}, &\forall P_{ij} \in \cP - \cC, \, 0 \leq n \leq N  \label{eq:PB1}\\
&p_{iju}^{0n} = p_{iju}^{Mn} &\forall P_{ij} \in \cC,\, 0 \leq n \leq k \label{eq:PB2}\\
&p_{ijl}^{0n} = p_{ijl}^{Mn} & \forall P_{ij} \in \cC,\, k \leq n \leq N \label{eq:PB3}\\
&R^0_{ij} = R^M_{ij}& \forall P_{ij} \in \cC. \label{eq:PB4}
\end{align}
The integral in the objective of problem \eqref{prob:ocp} is approximated by a Riemann sum (normalized by $T$) of the form
\begin{align} \label{discobj}
C_1\approx\sum_{P_{ij}\in\cC} \sum_{m=0}^M \frac{1}{M+1}S_{ij}^m.
\end{align}
Finally, we add to the objective a penalty on total squared variation of compressor ratio values, in order to guarantee a physically realistic solution of the nonlinear program.  This penalty takes the form
\begin{align}
C_2\approx \sum_{P_{ij}\in\cC} \sum_{m=0}^{M-1} & \frac{1}{M} (R_{ij}^m - R_{ij}^{m+1})^2. \label{obj:smoothTrap}
\end{align}
The nonlinear program is then given by
\begin{equation} \label{prob:ocp2}
\begin{array}{llll}
\!\!\!\!\!\!\! \min & \dS C_1+\mu C_2 \text{ using } \eqref{discobj}, \eqref{obj:smoothTrap} \\ \!\!\!\!\!\!\! s.t.
& \text{pipe dynamics: }  \eqref{eq:Trape1}-\eqref{eq:Dis2} \,\, \& \,\, \eqref{eq:Trape4}-\eqref{eq:Dis6}  \\
& \text{joint conditions: }  \eqref{bal5},\eqref{bal6}\\
& \text{density \& compression constraints: } \eqref{eq:PBound1}-\eqref{compdisc} \\
& \text{periodicity constraints: } \eqref{eq:PB1}-\eqref{eq:PB4} \\
& \text{injection parameters: } \eqref{bal7},\eqref{bal8} \\
& \text{compressor power: } \eqref{obj2}-\eqref{poscost}
\end{array}
\end{equation}
where $\mu$ is an appropriate scaling factor.  The software implementation of this nonlinear programming approach is described in the next section, followed by descriptions of computational experiments involving several examples with comparisons to the results of other studies.

\section{Implementation and Examples} \label{sec:results}

Instead of solving the above nonlinear program directly, the
implementation employs a lexicographic strategy. It first solves the
nonlinear program with the original objective \eqref{discobj}.  It
then solves the nonlinear program with the smoothness objective
\eqref{obj:smoothTrap}, while imposing the additional constraint
\begin{align}
 C_1 \leq (1+r) f, \mbox{ where } 0 \leq r \leq 1
\end{align}
where $f$ is the optimal objective value of the first
step. Intuitively, the tolerance $r$ is a user-adjustable parameter
that quantifies the percentage of compression energy that can be traded for
a smoother solution. Note also that the second optimization stage is initialized
using the first-stage solution.

This large-scale nonlinear problem is modeled with AMPL (version 2014)
\cite{fourer02ampl,amplLink} and solved with the nonlinear solver
IPOPT 3.12.2, ASL routine (version 2015) \cite{wachter06on}.  The
implementation is run on a Dell PowerEdge R415 with AMD Opteron 4226
and 64 GB ram.

We present the computational results of three case studies that include a validation
of the proposed approach, as well as results about its solution quality, efficiency, and scalability.


\begin{figure}[!t]
\centering
\includegraphics[width=.9\linewidth]{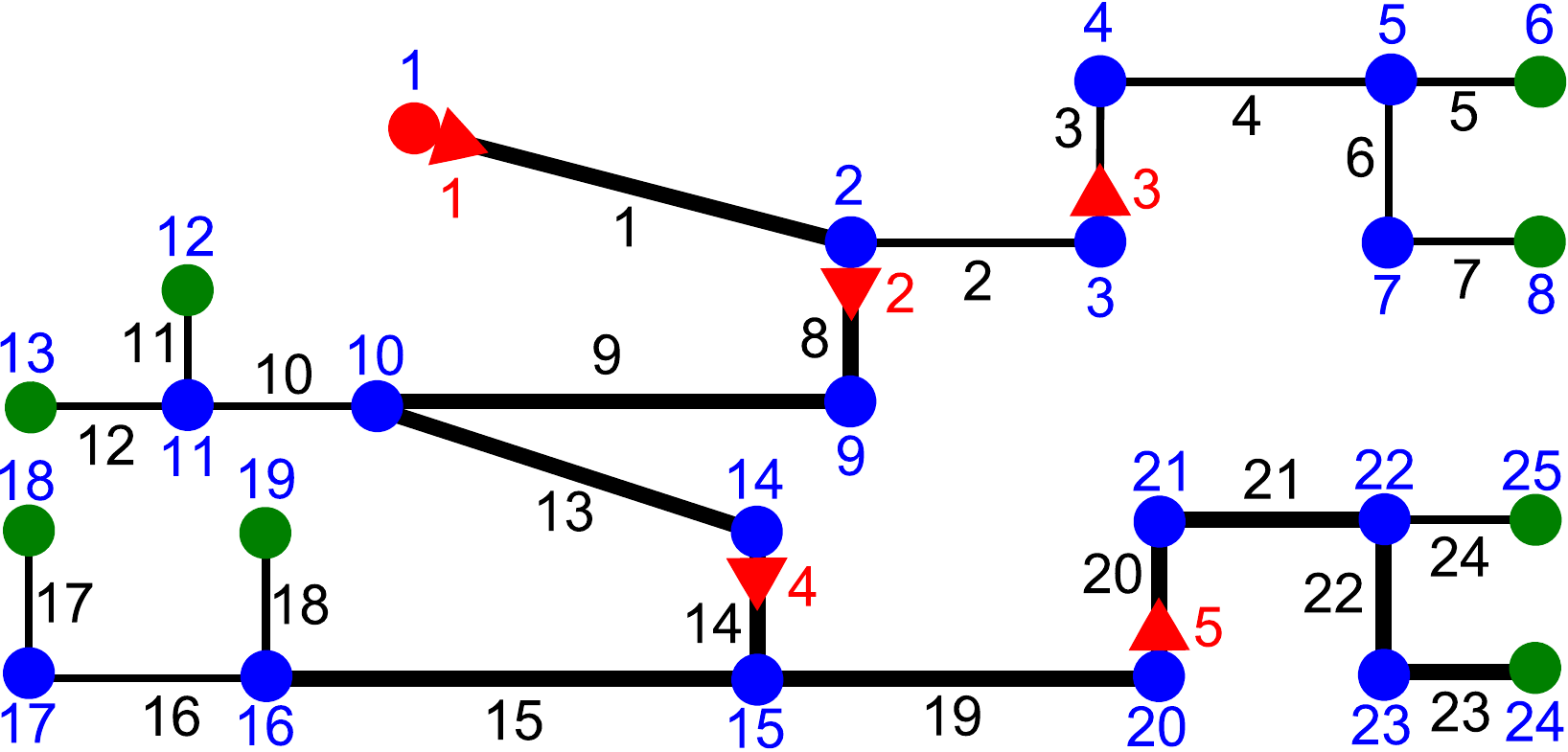}
\caption{24-pipe gas system test network used in the benchmark case study.  Numbers indicate nodes (blue), edges (black), and compressors (red). Thick and thin lines indicate 36 and 25 inch pipes.  Nodes are source (red), transit (blue), and consumers (green).}
\label{fig:model9} \vspace{-3ex}
\end{figure}

\paragraph{Validation}

The solution obtained using our implementation was validated on the 24-pipe benchmark gas network used in prior work \cite{zlotnik15cdc}, and illustrated in Figure \ref{fig:model9}. The pressures at supply sources were fixed at 500psi, the scalings for non-dimensionalization were set to $p_N$ = 250psi and $q_N$ = 100 kg/m$^2$/s, physical parameters $a$ =
377.968 m/s, $\gamma = 2.5$, and $\lambda = 0.01$ were used, and a time horizon $T = 86400$ seconds (24 hours) was considered.  Parameters $D_{ij}, A_{ij}, L_{ij}, \overline{R_{ij}},$ and $\underline{R_{ij}}$ were set according to the benchmark case study, as well as time-dependent profiles of gas injections/withdrawals $d_i(t)$.  The benchmark network structure and the gas draw profiles are provided online \cite{mak16data}.  Each pipe $P_{ij}$ is discretized uniformly according to its length $L_{ij}$ into $\lceil
L_{ij}/E \rceil + 1$ segments, where $E$ is set to 10km by default. In order to correspond to the structure of the benchmark model, the compressors are placed on the first segment of the $i^{\mbox{th}}$ end of every pipe $P_{ij} \in \cC$.

The admissible pressure range is 500 to 800 psi throughout the network.  A feasible solution to the discretized problem that satisfies the pressure constraints may cause these constraints to be violated in a high-accuracy simulation of the dynamics for the continuous problem.  To address this issue, one version of our implementation tightens the pressure bounds conservatively by 4\% or less, i.e. 520 to 780 psi for the benchmark case.  We refer to this as ``tightened'', and optimization using the nominal constraints of 500 to 800 psi is referred to as ``regular''.

\begin{figure*}[!t]
\centering
\includegraphics[width=.24\textwidth,height=2.7cm]{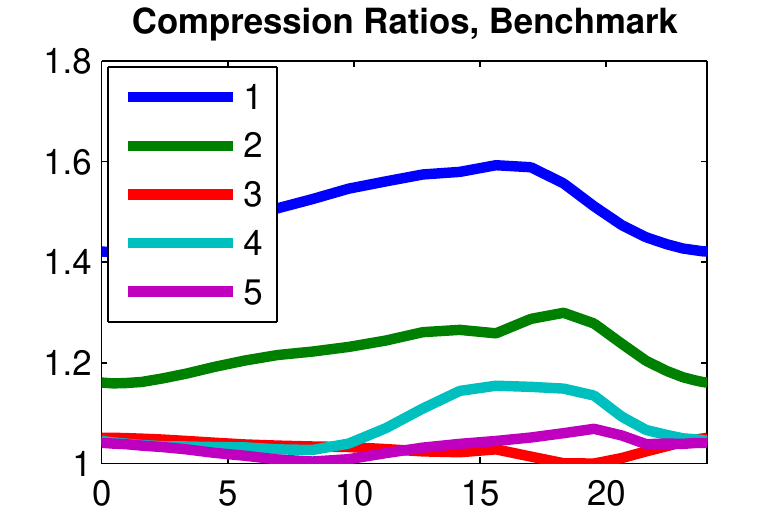}
\includegraphics[width=.24\textwidth,height=2.7cm]{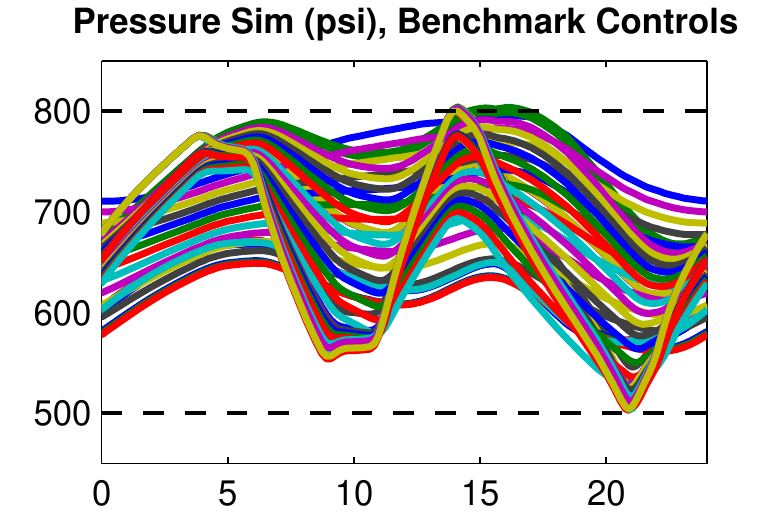}
\includegraphics[width=.24\textwidth,height=2.7cm]{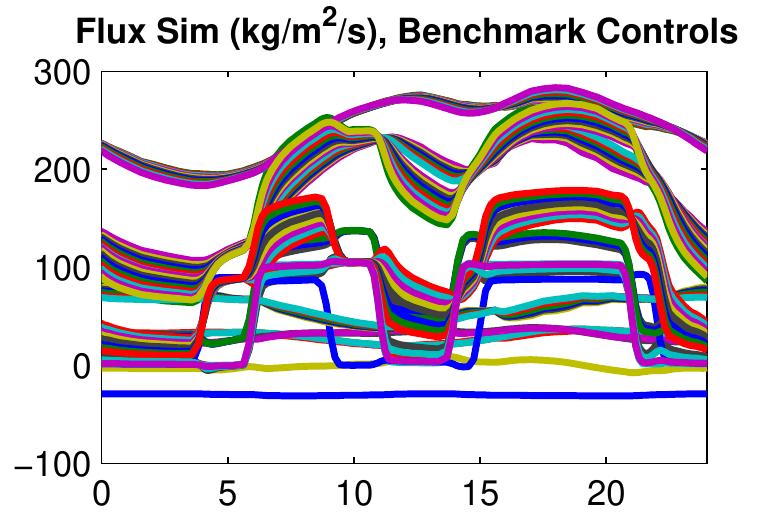}
\includegraphics[width=.24\textwidth,height=2.7cm]{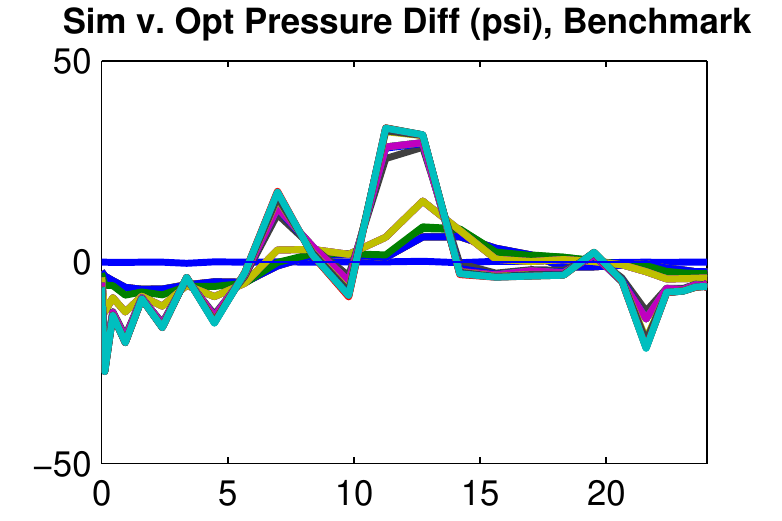} \\
\includegraphics[width=.24\textwidth,height=2.7cm]{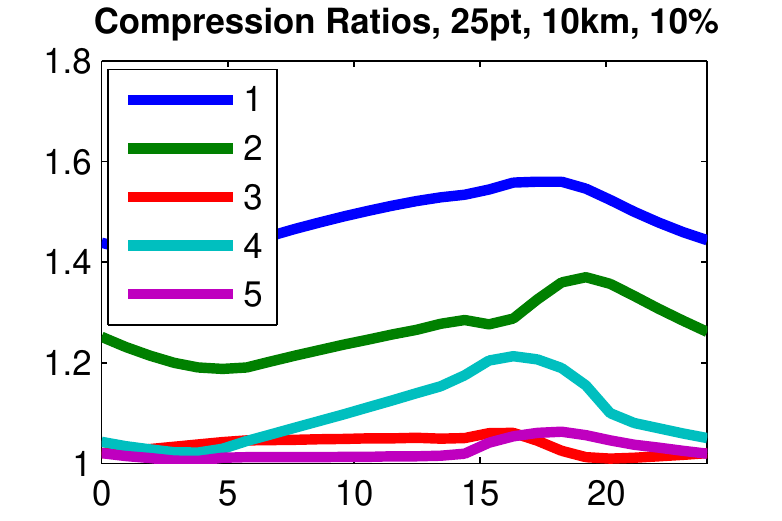}
\includegraphics[width=.24\textwidth,height=2.7cm]{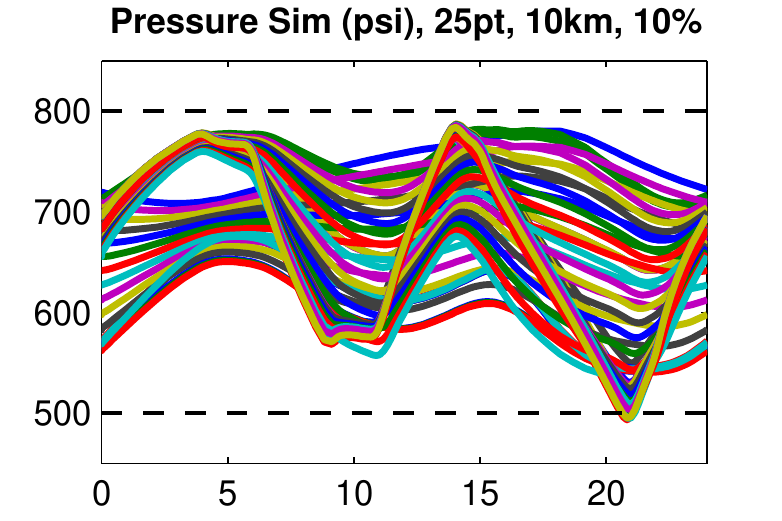}
\includegraphics[width=.24\textwidth,height=2.7cm]{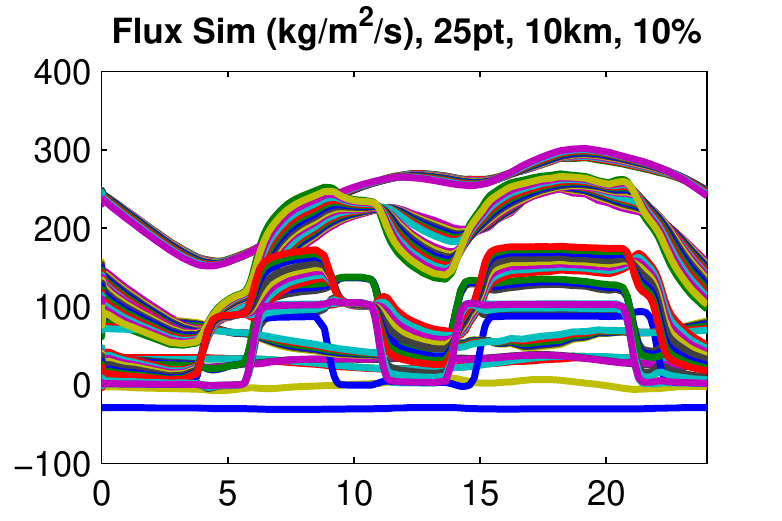}
\includegraphics[width=.24\textwidth,height=2.7cm]{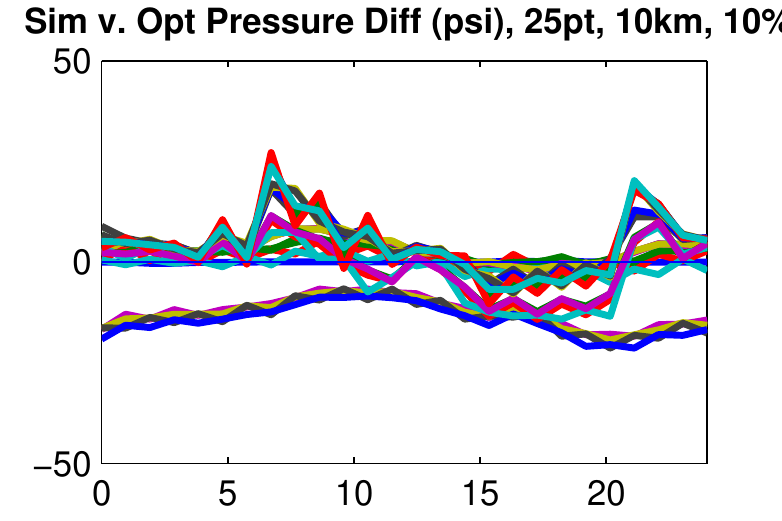} \\
\includegraphics[width=.24\textwidth,height=2.7cm]{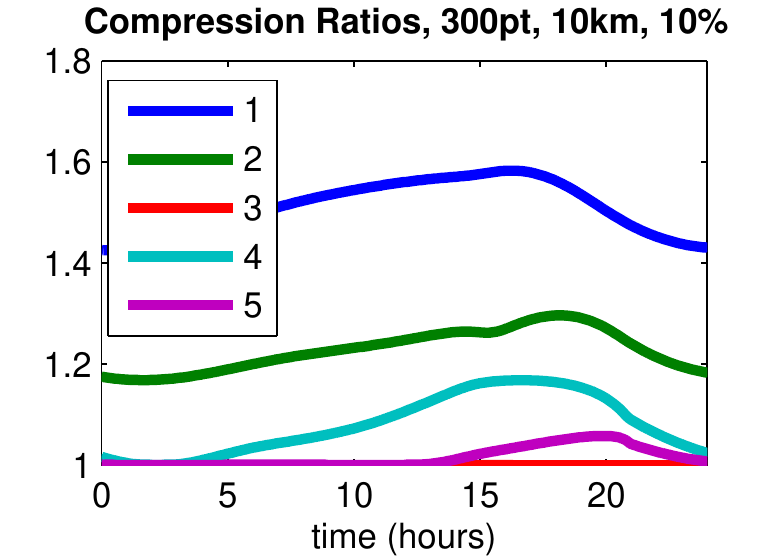}
\includegraphics[width=.24\textwidth,height=2.7cm]{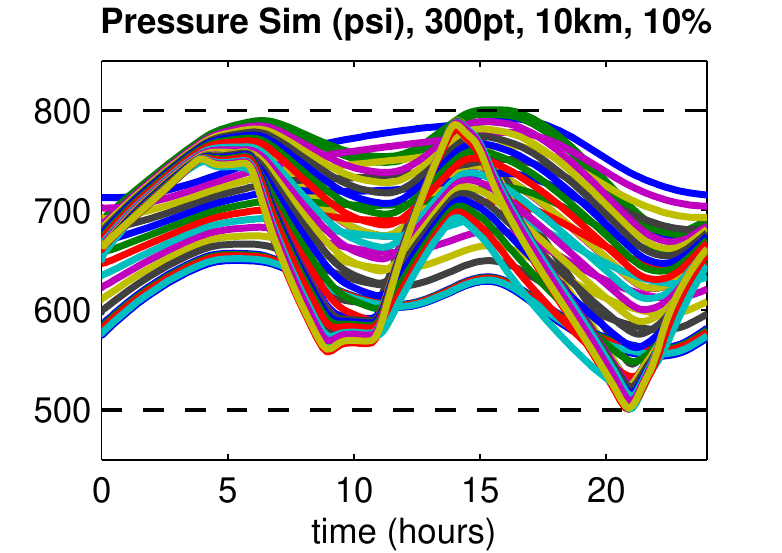}
\includegraphics[width=.24\textwidth,height=2.7cm]{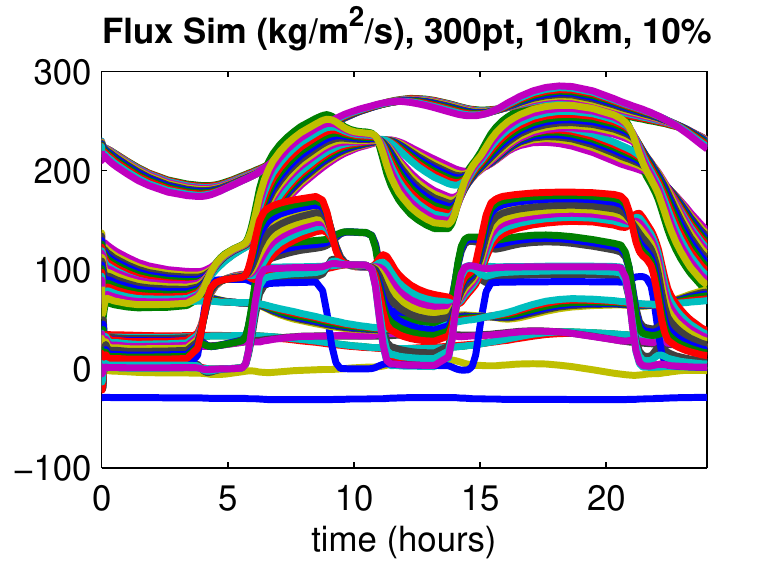}
\includegraphics[width=.24\textwidth,height=2.7cm]{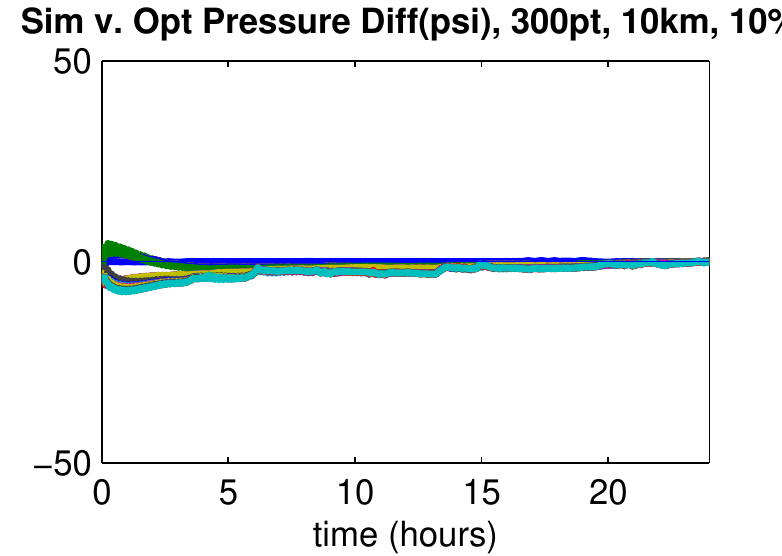} \\
\caption{From top to bottom: Benchmark solution; New method with 25 points in time, tightened constraints, 10km in space, $r=10\%$; New method with 300 points in time, regular constraints, 10km in space, $r=10\%$.  From left to right: Optimal control solution; Pressure trajectories from simulation using the controls; Flux trajectories from the same simulation; Difference between pressure solution from optimization and pressure trajectories from simulation.} \vspace{-2ex}
\label{fig:validate}
\end{figure*}

The optimization results were validated by using the optimized compression ratio solution as a time-varying parameter in a validated dynamic simulation method \cite{grundel13b,zlotnik15dscc}, and through comparison to a state-of-the-art (benchmark) solution based on a pseudospectral discretization approach \cite{zlotnik15cdc}.  The trajectories computed using the simulation are used to validate the optimization solution in two ways.  First, we quantify how much the constraints on pressure are exceeded by evaluating the $L_2$ norm of the violation.  This aggregates violations over the 24-hour period by integrating the square of pressure violation (psi) of the bounds at every junction. This takes the form
\begin{align} \label{ccompare1}
\!\!v_p=\sqrt{\sum_{P_{ij}\in \cP} V_{ij}^p}
\end{align}
where
\begin{align} \label{pviolnorm}
V_{ij}^p & = \bq{\int_0^T((p_{ij}(t,0)-p_{\max})_+)^2\rd t}^{\thalf} \nonumber \\ & \quad +\bq{\int_0^T((p_{\min}-p_{ij}(t,L_{ij}))_+)^2\rd t}^{\thalf}
\end{align}
and $(x)_+=x$ if $x\geq 0$ and $(x)_+\equiv 0$ if $x<0$.  The unit of the metric is psi-days.  Table \ref{sim:violations} lists its values for solutions obtained using various time discretization $M$, smoothing parameter $r$, and using tightened and regular constraints.  With the tightened bounds, the optimization solution has no, or very small violations for these configurations. Figure \ref{fig:validate} depicts the optimal compressor ratio functions for 25 and 300 point temporal discretization, with tightened and regular constraints, respectively, and both with $E=10$km spatial discretization and $r=10\%$.  A comparison is made with the solution to the same network obtained from prior work \cite{zlotnik15cdc}.  The smoothness of the compressor ratios over time indicates that the optimization solution is physically realistic and the control profiles can be implemented by operators. These results suggest that our optimization method provides operationally feasible solutions.

Next, we further validate the solutions provided by the solver by comparing the optimized pressure profiles with pressure trajectories obtained by applying the optimal controls in a dynamic simulation of the ODE model of the network \cite{zlotnik15cdc} using an adaptive time-stepping solver \texttt{ode15s} in MATLAB.  Figure \ref{fig:validate} depicts the pressure and flow profiles resulting from the simulation, and the difference between optimization and simulation pressure trajectories is illustrated.  For the 25 point time discretization, the discrepancy is similar to that for the benchmark solution, and using a 300 point time discretization nearly eliminates this difference.  Because the sources of the compared pressure profiles are qualitatively very different, i.e., optimization of algebraic equations that discretize PDEs over a fixed grid compared with adaptive time-stepping solution of an ODE system, this is a powerful cross-validation of both models.

Finally, Table \ref{sim:24B} presents the aggregated differences (in relative L2-Norm) between compressor ratios of the proposed method and the
state-of-the-art benchmark solution. Specifically, we compute
\begin{align} \label{ccompare1}
\!\!v_R=\sqrt{\sum_{P_{ij}\in \cC} (V_{ij}^R)^2}
\end{align}
where
\begin{align} \label{ccompare2}
\dS V_{ij}^{R}  = \frac{\int_0^T(R_{ij}(t)-R_{ij}^*(t))^2\rd t}{\int_0^T(R_{ij}^*(t))^2\rd t}.
\end{align}
Once again, the results indicate that the proposed method and the state-of-the-art yield highly similar compressor ratios across the network.

\begin{table}[t]
  \centering
  \caption{Aggregated Pressure Bound Violations ($v_p$, psi-days).}
  \scriptsize
    \begin{tabular}{c|cc|cc}
    \toprule
    Bounds & \multicolumn{2}{c}{Tightened} & \multicolumn{2}{c}{Regular}  \\
    \midrule
    Time pt. & 5\% & 10\% & 5\% & 10\% \\
    \midrule
      TZ 25tp   &0.0922 &0.1296 & 1.4994 & 1.9574  \\
      TZ 50tp   & 0.0000 &0.0000 & 0.2987 & 0.2693 \\
      TZ 100tp &0.0000 &0.0000 & 0.0394 & 1.0681 \\
      TZ 300tp     &0.0000 &0.0000 & 0.0062 & 0.0414 \\
    \midrule
  \end{tabular}
  \label{sim:violations}%
\end{table}

\begin{table}[t]
  \centering
  \scriptsize
  \caption{Distance ($v_R$) of Compression Ratios w.r.t. Benchmarks}
    \begin{tabular}{c|cc|cc}
    \toprule
    Bounds & \multicolumn{2}{c}{Tightened} & \multicolumn{2}{c}{Regular} \\
    \midrule
    Time pt. &  5\% & 10\%  &  5\% & 10\% \\
    \midrule
      TZ 25tp   &0.0013& 0.0015 & 0.0023 & 0.0015 \\
      TZ 50tp &0.0018  &0.0027 & 0.0022 & 0.0011  \\
      TZ 100tp & 0.0018 & 0.0031 & 0.0018 & 0.0014 \\
      TZ 300tp  & 0.0006 &   0.0006 & 0.0009 & 0.0007\\
      \bottomrule
    \end{tabular}%
  \label{sim:24B}%
\end{table}%


%


\paragraph{Solution Quality and Efficiency}

Table \ref{2stage:24} reports the objective value $C_1$ and computation time
of the proposed method for various time discretizations and smoothness
parameters $r$.  The number of variables in the optimization problem ranges from 11,675 for the 25pt discretization to 135,150 for the 300pt model.  The table gives the value of the $C_1$ objective after
the first stage, and also in the second stage for $r=5\%$ and
$10\%$.  The CPU times in seconds for the first and
second stages are also reported. First, we observe that enforcing the smoothness of the solution does not fundamentally decrease the quality of the $C_1$ objective, which is important from an operational standpoint. Second, as expected, refining the time discretization increases the objective value, because more constraints are added.  The
convergence rate is fast and the solutions obtained with a coarse
discretization are already of high quality, as illustrated in Figure \ref{fig:validate}. Finally, the method is extremely efficient; conisider the time granularities with 25 and 40 points: For $r=10\%$, the method takes about 20 and 80 seconds, respectively, which indicates that it can be used during real-time
operations.  The line (PS) reports the efficiency results of our
proposed method when a pseudospectral (PS) discretization
\cite{elnagar1995pseudospectral,gong2008spectral} is used in time instead of
the trapezoidal discretization.  This computation is conducted for calibration purposes to show that the proposed method is at least one order of magnitude
faster than a pseudospectral discretization method when a two-stage strategy is used to ensure a smooth, physically realistic solution. This is not
surprising, because the resulting optimization problem has a much denser
matrix and requires more iterations, thus decreasing the overall solver efficiency.  The last line describes the objective value and the execution of the state-of the art (benchmark) solution \cite{zlotnik15cdc}, which was executed on a different machine, and only given to indicate that our implementation of the pseudospectral discretization is comparable to previous work.

\begin{table}[t]
  \centering
  \scriptsize
  \caption{Objective Value ($C_1$) and Runtimes on 24 Pipe Network.}
    \begin{tabular}{c|c|cc|c|cc}
    \toprule
             & \multicolumn{3}{c|}{Objective Value} & \multicolumn{3}{c}{CPU Time (secs)} \\
    \midrule
     & 1st Stage & \multicolumn{2}{c|}{2nd Stage} & 1st Stage & \multicolumn{2}{c}{2nd Stage} \\
    \midrule
    &            & $r= $ 5\%   & 10\%  &   &  $r= $ 5\%   & 10\%  \\
    \midrule
    TZ 25tp &4.088 & 4.292 & 4.496  & 8 & 25 & 14  \\
    TZ 40tp &4.235 & 4.447&  4.658  &41 & 95 &37 \\
    TZ 50tp &4.231 & 4.443&  4.654  & 32 & 122 & 49 \\
    TZ 60tp &4.304 & 4.519 & 4.734  &79 &  104&78  \\
    TZ 80tp &4.330 & 4.546&  4.763  & 79& 153 & 81 \\
    TZ 100tp &4.332 & 4.548 &4.765  &135 & 290 & 139 \\
    TZ 300tp &4.333 & 4.550 &4.766  & 1041& 452 & 435\\
    \midrule
    PS 25tp &       &       &       & 252 & 266 &394 \\
    \midrule
    Benchmark & 4.091 &       &       & 253$^*$ &  & \\
    \bottomrule
    \end{tabular}%
  \label{2stage:24}%
\end{table}%

\paragraph{Scalability}

\begin{table}[t]
  \centering
  \scriptsize
  \caption{Objective Value ($C_1$) and Runtimes on 40 Pipe Network.}
    \begin{tabular}{c|c|cc|c|cc}
    \toprule
             & \multicolumn{3}{c|}{Objective Value} & \multicolumn{3}{c}{CPU Time (secs)} \\
    \midrule
     & 1st Stage & \multicolumn{2}{c|}{2nd Stage} & 1st Stage & \multicolumn{2}{c}{2nd Stage} \\
    \midrule
    &            & $r= $ 10\%   & 20\%  &   &  $r= $ 10\%   & 20\%  \\
    \midrule
    TZ 20pt & 0.528 &  0.580 & 0.633 &  69  &  20  & 43\\
    TZ 25tp & 0.558 &  0.614 & 0.670 &  220 &  45  & 48\\
    TZ 40tp & 0.596 &  0.655 & 0.715 &  846 & 305  & 757\\
    \midrule
    PS 20tp &       &        &       & 372  & 376  &577\\
    \bottomrule
    \end{tabular}%
  \label{2stage:40}%
\end{table}%

\begin{table}[t]
  \centering
  \scriptsize
  \caption{Objective Value ($C_1$) and Runtimes on 135 Pipe Network.}
    \begin{tabular}{c|c|cc|c|cc}
    \toprule
             & \multicolumn{3}{c|}{Objective Value} & \multicolumn{3}{c}{CPU Time (secs)} \\
    \midrule
     & 1st Stage & \multicolumn{2}{c|}{2nd Stage} & 1st Stage & \multicolumn{2}{c}{2nd Stage} \\
    \midrule
    &            & $r= $ 5\%   & 10\%  &   &  $r= $ 5\%   & 10\%  \\
    \midrule
    TZ 10tp & 3.406&  3.577  &  3.747 & 2543   & 1507  &  910 \\
    TZ 15tp & 4.002&  4.202  &  4.402 & 4415   & 1313  & 1674  \\
    TZ 20tp & 4.735 & 4.971  &  5.208 & 16674  & 11289 &  2858 \\
    TZ 25tp & 4.295 & 4.510  &  4.725 & 37143  &2270   & 5249  \\
    \bottomrule
    \end{tabular}%
  \label{2stage:135}%
\end{table}%

\begin{figure}[t]
\center
\includegraphics[width=.49\linewidth]{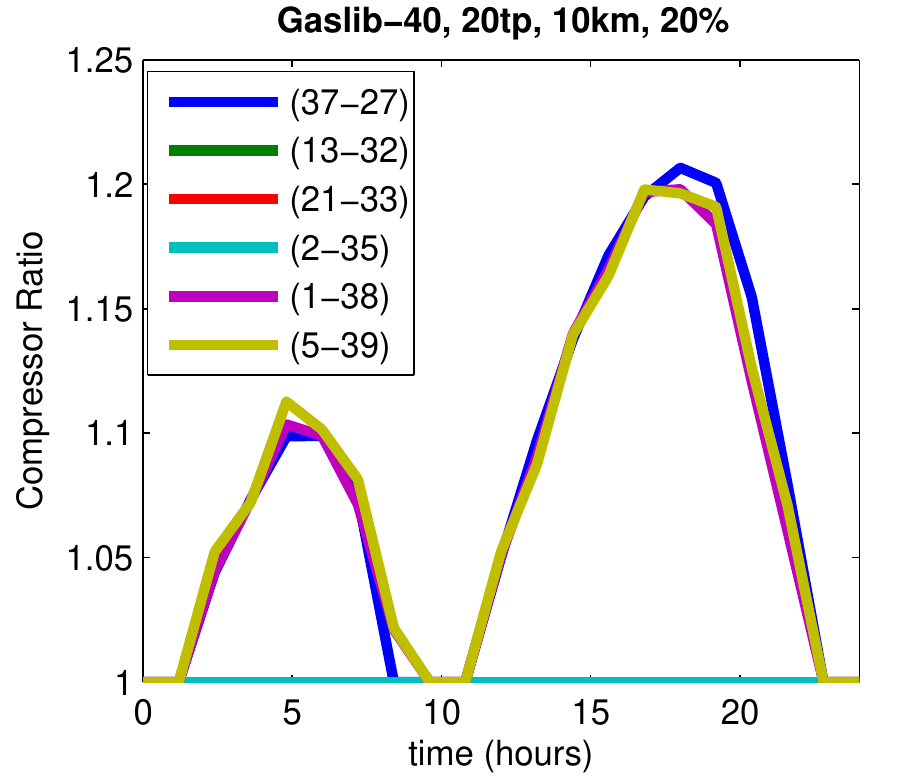}
\includegraphics[width=.49\linewidth]{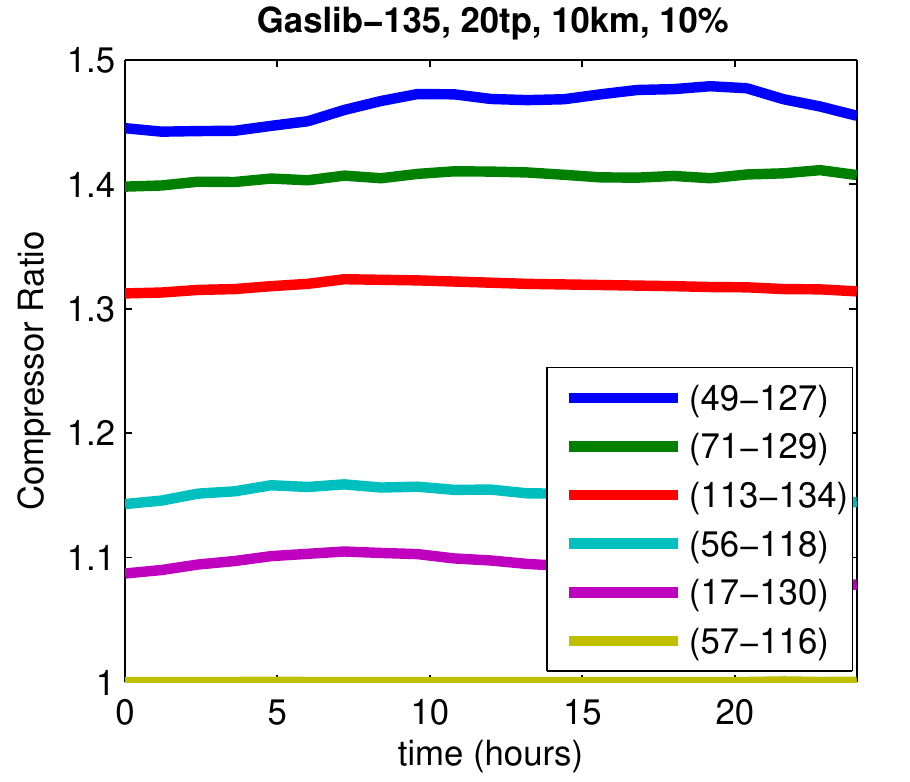}
\caption{Compression ratio solutions for gaslib-40 (left) and gaslib-135 (right) case studies.  The trapezoidal scheme is used for time discretization with 20 points, and the smoothness tolerances are $r=20\%$ and $r=10\%$, respectively.} \vspace{-3ex}
\label{fig:comprats}
\end{figure}

To study how the proposed method scales to larger network models, two
additional instances are considered: gaslib-40 and gaslib-135 from the
GasLib library \cite{pfetsch2015validation}.  The pressure ranges are
set to 500 to 800 psi and 500 to 1000 psi for gaslib-40 and gaslib-135,
respectively, and the source pressure is set to 600 psi.  Simplified network structures and time-dependent gas draw profiles are provided online \cite{mak16data}.  Tables
\ref{2stage:40} and \ref{2stage:135} present the results, which
show that the proposed method scales well to larger benchmarks and
exhibits similar behavior as on the 24-pipe network. The 40 pipe
network is solved in less than two minutes, and the 135 pipe network in
less than an hour. The pseudo-spectral discretization does not
converge when applied to the 135 pipe network, which explains why no
results were presented. Figure \ref{fig:comprats} presents the
smoothness results for these two case studies. Benchmark gaslib-40 is
more challenging in this respect, so that higher tolerances
are required. Only a subset of the compressors are shown for these cases.

\section{Conclusions}
\label{sec:conc}

We have investigated the Dynamic Optimal Gas Flow (DOGF) problem in which
the goal is to minimize compressor costs under dynamic conditions
where deliveries to customers are described by time-dependent mass
flow functions.  This study is motivated by the growing reliance of electric power systems on gas-fired generation in order to balance intermittent
sources of renewable energy.  This deeper integration has increased
the variation and volume of flows through natural gas transmission
pipelines, emphasizing the need to go beyond steady-state optimization
in controlling and optimizing natural gas networks. Maintaining
efficiency and security under such dynamic conditions requires
optimization methods that account for transients and can quickly
compute solutions to follow generator re-dispatch.

This paper presents an efficient scheme for the DOGF that relies on a
compact representation of gas flow physics and a trapezoidal
discretization in time and space that allows sparse representation of constraints.  A two-stage approach is applied to minimize energy costs and maximize smoothness  to correspond to physically realistic solutions. The resulting large-scale
nonlinear programs are solved using a modern interior-point method and
the results are validated against an accurate simulation
of the dynamic equations and a recently proposed state-of-the-art
method. The novel optimization scheme yields solutions that are
feasible for the continuous problem and practical from an operational
standpoint in a computation time that is at least an order of magnitude faster
than existing methods. Scalability of the scheme is
demonstrated using three networks with 24, 40, and 135 pipes with total lengths of 477, 1118, and 6964 km of pipeline, respectively. These results indicate that optimization of dynamic flows through gas pipeline networks is now in the realm of optimization technology.  This opens new directions for the investigation of previously intractable control and optimization problems involving such networks, as well as integration of their operations with electric power systems.

Future work will further accelerate computation by
exploiting the underlying time structure of the optimization problem
and considering alternative discretization schemes, and how
they impact accuracy and performance.


{\small
\bibliographystyle{unsrt}
\bibliography{gas_master,power_master,TS}
}

\end{document}